\documentclass{amsart}

\theoremstyle{plain}
\newtheorem{theorem}{Theorem}[section]
\newtheorem{lemma}{Lemma}[section]
\newtheorem{proposition}{Proposition}[section]
\numberwithin{equation}{section}
\newcommand{\Z}{\mathbb{Z}}
\newcommand{\qbin}[2]{\genfrac{[}{]}{0pt}{}{#1}{#2}}
\newcommand{\qbins}[2]{{\textstyle\genfrac{[}{]}{0pt}{}{#1}{#2}}}

\begin{document}

\title[Rogers--Ramanujan identities]{Partial-sum analogues of the
Rogers--Ramanujan identities}

\author[Ole Warnaar]{S. Ole Warnaar}

\address{Department of Mathematics and Statistics,
The University of Melbourne, Vic 3010, Australia}
\email{warnaar@ms.unimelb.edu.au}


\thanks{Work supported by the Australian Research Council}

\dedicatory{Dedicated to Barry McCoy on the occasion of his sixtieth birthday}

\begin{abstract}
A new polynomial analogue of the Rogers--Ramanu\-jan identities 
is proven. Here the product-side of the Rogers--Ramanu\-jan identities is
replaced by a partial theta sum and the sum-side by a weighted sum over
Schur polynomials.
\end{abstract}


\maketitle

\section{Introduction}
The famous Rogers--Ramanujan identities are given by \cite{Rogers94}
\begin{equation}\label{RR1}
1+\sum_{n=1}^{\infty} \frac{q^{n^2}}{(1-q)(1-q^2)\cdots(1-q^n)}=
\prod_{j=0}^{\infty}\frac{1}{(1-q^{5j+1})(1-q^{5j+4})}
\end{equation}
and
\begin{equation}\label{RR2}
1+\sum_{n=1}^{\infty} \frac{q^{n(n+1)}}{(1-q)(1-q^2)\cdots(1-q^n)}
=\prod_{j=0}^{\infty}\frac{1}{(1-q^{5j+2})(1-q^{5j+3})}
\end{equation}
for $|q|<1$.
In one of his two proofs of these identities, Schur \cite{Schur17}
introduced two sequences of polynomials $(e_n)_{n\geq 2}$ and 
$(d_n)_{n\geq 2}$, where $e_n$ ($d_n$) is the generating function of 
partitions with difference between parts at least $2$ (and no part
equal to $1$), and largest part at most $n-2$.
The partitions $\{\emptyset,(1),(2),(3),(4),(3,1),(4,1),(4,2)\}$, for example, 
contribute to $e_6=1+q+q^2+q^3+2q^4+q^5+q^6$ and
the partitions $\{\emptyset,(2),(3),(4),(4,2)\}$ contribute to
$d_6=1+q^2+q^3+q^4+q^6$.

By standard combinatorial arguments, see e.g., \cite{MacMahon16,GR90},
it follows that $e_{\infty}:=\lim_{n\to\infty} e_n
=\text{LHS}\eqref{RR1}$
and $d_{\infty}:=\lim_{n\to\infty} d_n=\text{LHS}\eqref{RR2}$. 
Schur proved the Rogers--Ramanujan identities by showing that these
limits also hold when $\text{LHS}$ is replaced by $\text{RHS}$.
This he achieved by showing that both $e_n$ and $d_n$ satisfy the recurrence
\begin{equation}\label{xrec}
x_{n+2}=x_{n+1}+q^n x_n,
\end{equation}
and by solving this recurrence subject to the initial conditions $d_1=0$, 
$e_1=e_2=d_2=1$ (consistent with the 
combinatorial definition of $e_n$ and $d_n$ for $n\geq 2$). 
Specifically, Schur's solution to \eqref{xrec} reads
\begin{subequations}\label{edbos}
\begin{align}
e_n&=\sum_{j=-\infty}^{\infty}(-1)^j q^{j(5j+1)/2}
\qbin{n-1}{\lfloor{(n-5j-1)/2}\rfloor} \\
d_n&=\sum_{j=-\infty}^{\infty}(-1)^j q^{j(5j+3)/2}
\qbin{n-1}{\lfloor{(n-5j-2)/2}\rfloor}
\end{align}
\end{subequations}
for $n\geq 1$ and $\lfloor x \rfloor$ denoting the integer part of $x$.
Here the $q$-binomial coefficients are given by 
$\qbins{n}{m}=(q;q)_n/(q;q)_m(q;q)_{n-m}$ for 
$0\leq m\leq n$ and zero otherwise,
where $(a;q)_n=\prod_{j=0}^{n-1}(1-aq^j)$.

Employing the notation $(a_1,\dots,a_k;q)_n=
(a_1;q)_n\cdots(a_k;q)_n$ and recalling to the Jacobi triple product 
identity \cite[Eq. (II.28)]{GR90}
\begin{equation}\label{tpi}
\sum_{n=-\infty}^{\infty}(-1)^n a^n q^{\binom{n}{2}}=(a,q/a,q;q)_{\infty}
\end{equation}
it is now easy to obtain the desired limits;
\begin{align*}
e_{\infty}&=\frac{1}{(q;q)_{\infty}}
\sum_{j=-\infty}^{\infty}(-1)^j q^{j(5j+1)/2}=\frac{1}{(q,q^4;q^5)_{\infty}}
=\text{RHS}\eqref{RR1} \\
d_{\infty}&=\frac{1}{(q;q)_{\infty}}
\sum_{j=-\infty}^{\infty}(-1)^j q^{j(5j+3)/2}=\frac{1}{(q^2,q^3;q^5)_{\infty}}
=\text{RHS}\eqref{RR2}.
\end{align*}

Representations for the Schur polynomials similar to the left sides
of the Rogers--Ramanujan identities are also known 
\cite[\S 286 and \S 289]{MacMahon16},
\begin{equation}\label{edfer}
e_n=\sum_{r=0}^{\infty}q^{r^2}\qbin{n-r-1}{r} \quad\text{and}\quad
d_n=\sum_{r=0}^{\infty}q^{r(r+1)}\qbin{n-r-2}{r}.
\end{equation}
Equating this with \eqref{edbos} yields the following
polynomial analogue of the Rogers--Ramanujan identities \cite{Andrews70}:
\begin{equation}\label{RRpoly}
\sum_{r=0}^{\infty}q^{r(r+a)}\qbin{n-r-a}{r}= 
\sum_{j=-\infty}^{\infty}(-1)^j q^{j(5j+2a+1)/2}
\qbin{n}{\lfloor{(n-5j-a)/2}\rfloor}
\end{equation}
for $n\geq 0$ and $a\in\{0,1\}$.

Recently there has been renewed interest in the Schur polynomials
\cite{GIS99,AKP00,IPS00,BP01,W01} sparked by the following
nice generalization of the Rogers--Ramanujan identities due to
Garrett, Ismail and Stanton \cite[Eq. (3.5)]{GIS99}
\begin{equation}\label{GIS}
\sum_{r=0}^{\infty}\frac{q^{r(r+m)}}{(q;q)_r}=
\frac{(-1)^m q^{-\binom{m}{2}}d_m}{(q,q^4;q^5)_{\infty}}-
\frac{(-1)^m q^{-\binom{m}{2}}e_m}{(q^2,q^3;q^5)_{\infty}},
\end{equation}
where $m$ is a nonnegative integer and
$e_0=0$ and $d_0=1$ consistent with \eqref{xrec}. 

In this paper we show that \eqref{GIS} may be used to prove
new polynomial analogues of the Rogers--Ramanujan identities involving
the Schur polynomials. These polynomial identities are fundamentally
different from \eqref{RRpoly} in that the product-side is replaced by a 
partial theta series.
\begin{theorem}\label{thm}
For $k\in\{0,1\}$ and $n\geq 0$ there holds
\begin{align} \notag 
\sum_{j=-n-k}^n&(-1)^j q^{j(5j+1)/2} \\ 
&=\sum_{r=0}^n e_{2r+k+2}(-1)^{n-r}
q^{(n-r)(5n+3r+4k+5)/2}\frac{(q;q)_{n+r+k}}{(q;q)_{n-r}} 
\label{RRP1}\\ 
\intertext{and}
\notag 
\sum_{j=-n-k}^n&(-1)^j q^{j(5j+3)/2} \\
&=\sum_{r=0}^n d_{2r+k+2}(-1)^{n-r}
q^{(n-r)(5n+3r+4k+5)/2}\frac{(q;q)_{n+r+k}}{(q;q)_{n-r}}. 
\label{RRP2}
\end{align}
\end{theorem}
Partial theta-sum identities of this type were first discovered by 
Shanks \cite{Shanks51}.

When $n$ tends to infinity (for $|q|<1$) only the term with $r=n$ 
contributes to the sums on the right. Hence the first identity of the 
theorem implies
\begin{equation*}
\sum_{j=-\infty}^{\infty}(-1)^j q^{j(5j+1)/2}=
(q;q)_{\infty}e_{\infty}=(q;q)_{\infty}
\sum_{n=0}^{\infty} \frac{q^{n^2}}{(q;q)_n},
\end{equation*}
which is transformed into \eqref{RR1} by the triple product identity.
Likewise, \eqref{RR2} arises as the the large $n$ limit of the second
identity of the theorem.

Polynomial analogues of the Rogers--Ramanu\-jan strikingly similar to
those of Theorem~\ref{thm} have previously been discovered by 
Andrews \cite{Andrews83}.
For $n\geq 0$ let $K_n(x)$ denote the Szeg\"o polynomial~\cite{Szego26}
\begin{equation*}
K_n(x)=\sum_{r=0}^n x^r q^{r(r+1)}\qbin{n}{r}.
\end{equation*}
Then Andrews posed in the problems section of SIAM Review~\cite{Andrews83} 
the problem of showing that 
\begin{multline}\label{AJ}
\sum_{j=-n-k}^n(-1)^j q^{j(5j+2k+1)/2} \\
=\sum_{r=0}^n K_r(q^{2n-2r+k-1})(-1)^{n-r}
q^{(n-r)(5n-3r+4k+1)/2}\frac{(q;q)_{n+k}}{(q;q)_{n-r}}
\end{multline}
for $k=\{0,1\}$ and $n\geq 0$. Note here that the left side of
\eqref{AJ} coincides with the left side of \eqref{RRP1} (\eqref{RRP2})
when $k=0$ ($k=1$).

The remainder of this paper is divided in two parts with 
section~\ref{secp} containing a proof and section~\ref{secd}
a discussion of Theorem~\ref{thm}.
In the first part of this discussion we examine two simple proofs of 
\eqref{AJ} found by Jordan and Andrews and indicate our failure in 
generalizing these to a proof Theorem~\ref{thm}. 
The second part of our discussion focuses on some of the combinatorial aspects 
of Theorem~\ref{thm}.

\section{Proof of Theorem \ref{thm}}\label{secp}
\subsection{A more general identity}
Key to the proof of Theorem~\ref{thm} is the following proposition.
\begin{proposition}\label{prop2}
For $k\in\{0,1\}$ and $|a|,|q|<1$ there holds
\begin{multline}\label{key}
\sum_{n=0}^{\infty}\frac{a^{2n} q^{n(n+k)}}{(q;q)_n}=
\frac{(a;q)_{\infty}^2}{(q;q)_{\infty}^3}\sum_{j=1}^{\infty}
(q^{2j-k},q^{5+k-2j},q^5;q^5)_{\infty} \\
\times \sum_{r=0}^{\infty}\frac{(-1)^{j+r+1}
q^{\binom{j+r}{2}}(1-q^{2r+k+1})(aq^{-r};q)_r}{(a;q)_{r+k+1}}.
\end{multline}
\end{proposition}
It is perhaps not immediately clear that the sums on the right
converge, but inspection of the potentially problematic
terms shows that for $k\in\{0,1\}$ and $j\geq 1$,
\begin{multline*}
O\Bigl(q^{\binom{j+r}{2}}(q^{5+k-2j};q^5)_{\infty}(aq^{-r};q)_r\Bigr) \\
=\begin{cases}
q^{(j-1)(j+10r+4k+6)/10} & \text{$j\equiv 1,k+4 \pmod{5}$} \\
q^{(j-2)(j+10r+4k+7)/10+r+1} & \text{$j\equiv 2,k+3 \pmod{5}$,}
\end{cases}
\end{multline*}
which shows that both sums on the right converge and that their order
is irrelevant.

Before proving Proposition~\ref{prop2} we will show how it 
implies Theorem~\ref{thm}. Starting point is the observation that
\begin{multline}\label{inffin}
\sum_{j=1}^{\infty}(-1)^j q^{\binom{j+r}{2}}
\frac{(q^{2j-k},q^{5+k-2j},q^5;q^5)_{\infty}}{(q;q)_{\infty}} \\
=\sum_{i=1}^2 \frac{(-1)^{i+r}q^{-\binom{2r+k+2}{2}}}
{(q^{i+2k},q^{5-i-2k};q^5)_{\infty}}
\sum_{j=-r-k}^r (-1)^j q^{j(5j+2i+4k-5)/2}.
\end{multline}
To prove this we use that for $f_j$ such that
$f_j=0$ if $j\equiv 3k \pmod{5}$ there holds
\begin{equation*}
\sum_{j=1}^{\infty}f_j=\sum_{i=1}^2\sum_{j=0}^{\infty}
(f_{5j+i}+f_{5j+5-i+k}).
\end{equation*}
This, together with the simple to verify identities 
\begin{gather*}
(q^{m+5n},q^{5-5n-m};q^5)_{\infty}=
(q^m,q^{5-m};q^5)_{\infty}(-1)^n q^{-nm-5\binom{n}{2}} \\[2mm]
\frac{(q^{2i-k},q^{5+k-2i},q^5;q^5)_{\infty}}{(q;q)_{\infty}}=
\frac{1}{(q^{i+2k},q^{5-i-2k};q^5)_{\infty}}, \;
i,k+1\in\{1,2\} \\[2mm]
\binom{i+r}{2}-(r+1)(5r+2i+4k)/2=-\binom{2r+k+2}{2}, \;
i,k+1\in\{1,2\} \\[2mm]
\sum_{i=1}^2 (-1)^i=0
\end{gather*}
and the Jacobi triple product identity \eqref{tpi}, yields
\begin{align*}
&\text{LHS}\eqref{inffin} \\ 
&\; = \sum_{i=1}^2 \frac{(-1)^i q^{\binom{i+r}{2}}}
{(q^{i+2k},q^{5-i-2k};q^5)_{\infty}}
\biggl(\sum_{j=-\infty}^{-2r-k-2}\!+\sum_{j=0}^{\infty}\,\biggr)
(-1)^j q^{j(5j+10r+2i+4k+5)/2} \\
&\; =\sum_{i=1}^2 \frac{(-1)^{i+r} q^{\binom{i+r}{2}-(r+1)(5r+2i+4k)/2}}
{(q^{i+2k},q^{5-i-2k};q^5)_{\infty}} \\ 
&\qquad \times\biggl[\,\sum_{j=-r-k}^r(-1)^j q^{j(5j+2i+4k-5)/2}
-(q^{i+2k},q^{5-i-2k},q^5;q^5)_{\infty}\biggr] \\
&\; =\text{RHS}\eqref{inffin}.
\end{align*}
After substituting \eqref{inffin} in \eqref{key} we obtain
\begin{multline*}
\sum_{n=0}^{\infty}\frac{a^{2n} q^{n(n+k)}}{(q;q)_n} \\
=\frac{(a;q)_{\infty}^2}{(q;q)_{\infty}^2}\sum_{r=0}^{\infty} \sum_{i=1}^2
\frac{(-1)^{i+1}q^{-\binom{2r+k+2}{2}}(1-q^{2r+k+1})(aq^{-r};q)_r}
{(q^{i+2k},q^{5-i-2k};q^5)_{\infty}(a;q)_{r+k+1}} \\ 
\times \sum_{j=-r-k}^r(-1)^j q^{j(5j+2i+4k-5)/2}.
\end{multline*}
Here the reader is warned that the order of the sums over $r$ and $i$ must
be strictly adhered to. Indeed, our earlier considerations about 
convergence and the fact that \eqref{inffin} is true, guarantee the not 
so obvious fact that after summing over $i$ the sum over $r$ converges.

Our next step removes any further convergence issues as we now
specialize $a=q^{m+1}$ with $m$ a nonnegative integer. The sum over $r$
then terminates, yielding
\begin{multline}\label{am}
\sum_{n=0}^{\infty}\frac{q^{n(n+2m+k+2)}}{(q;q)_n}
=\sum_{i=1}^2 \frac{(-1)^{i+1}}{(q^{i+2k},q^{5-i-2k};q^5)_{\infty}}\\
\times \sum_{r=0}^m 
\frac{q^{-\binom{2r+k+2}{2}}(1-q^{2r+k+1})}{(q;q)_{m-r}(q;q)_{m+r+k+1}}
\sum_{j=-r-k}^r(-1)^j q^{j(5j+2i+4k-5)/2}.
\end{multline}
Rewriting the left-hand side using the Garrett--Ismail--Stanton identity 
\eqref{GIS} gives
\begin{equation*}
\frac{d_{2m+k+2}}{(q,q^4;q^5)_{\infty}}-
\frac{e_{2m+k+2}}{(q^2,q^3;q^5)_{\infty}}
=(-1)^k q^{\binom{2m+k+2}{2}}\text{RHS}\eqref{am}.
\end{equation*}
Multiplying both sides by $(q;q)_{2m+k+1}$ this is of the form
\begin{equation*}
\frac{P(q)}{(q^2,q^3;q^5)_{\infty}}=\frac{Q(q)}{(q,q^4;q^5)_{\infty}}
\end{equation*}
with $P(q)$ and $Q(q)$ polynomials. An identity of this type can only
be true if $P(q)=Q(q)=0$, and we infer
\begin{align*}
e_{2m+k+2}&=q^{\binom{2m+k+2}{2}}
\sum_{r=0}^m \frac{q^{-\binom{2r+k+2}{2}}(1-q^{2r+k+1})}
{(q;q)_{m-r}(q;q)_{m+r+k+1}}\sum_{j=-r-k}^r (-1)^j q^{j(5j+1)/2} \\
d_{2m+k+2}&=q^{\binom{2m+k+2}{2}}
\sum_{r=0}^m \frac{q^{-\binom{2r+k+2}{2}}(1-q^{2r+k+1})}
{(q;q)_{m-r}(q;q)_{m+r+k+1}}\sum_{j=-r-k}^r (-1)^j q^{j(5j+3)/2}
\end{align*}
for $m\geq 0$.
All that remains is to invert these new representations of the
Schur polynomials. This is easily done recalling the Bailey transform 
\cite{Andrews84}, which states that if
\begin{equation}\label{ba}
\beta_n=\sum_{r=0}^n\frac{\alpha_r}{(q;q)_{n-r}(aq;q)_{n+r}}
\end{equation}
then
\begin{equation}\label{ab}
\alpha_n=(1-aq^{2n})\sum_{r=0}^n\frac{(-1)^{n-r}q^{\binom{n-r}{2}}
(aq;q)_{n+r-1}}{(q;q)_{n-r}}\,\beta_r.
\end{equation}
For later reference we remark that a pair of sequences $(\alpha,\beta)$
that satisfies \eqref{ba} (or, equivalently, \eqref{ab}) is called a Bailey
pair relative to $a$.

Since our expressions for the Schur polynomials take the
form \eqref{ba} with $a=q^{k+1}$, we may invoke \eqref{ab}
to find the identity claimed in Theorem~\ref{thm}.

\subsection{Proof of Proposition \ref{prop2}}
Our proof relies on the following lemma.
\begin{lemma}\label{lem1}
For $k\in\{0,1\}$ and $M$ and $n$ integers there holds
\begin{multline*}
\frac{q^{n(n+2)}}{(q;q)_{\infty}^3}
\sum_{j=1}^{\infty}\sum_{r=0}^{\infty}\sum_{l=0}^M (-1)^{M+j+r+1}
q^{\binom{j+r}{2}+\binom{M-l}{2}+\binom{l}{2}+l(r+k+1)} \\
\times q^{-r(M-l)-n(2j-k)} 
\frac{(1-q^{(2j-k)(2n+1)})(1-q^{2r+k+1})}
{(1-q^{2n+1})(q;q)_{M-l}(q;q)_l}  \\[2mm]
=\begin{cases}\displaystyle
\frac{q^{m(m+k)}}{(q;q)_{m-n}(q;q)_{m+n+1}} & \text{$M=2m$} \\[4mm]
0 &\text{$M=2m+1$.} 
\end{cases}
\end{multline*}
\end{lemma}
Here our earlier definition of $(a;q)_n$ is extended to all integers $n$ by
$(a;q)_n=(a;q)_{\infty}/(aq^n;q)_{\infty}$. Note in particular that
$1/(q;q)_n=0$ for $n<0$.

Given the triple sum on the left, Lemma~\ref{lem1} perhaps appears 
complicated and not readily applicable. However, in view of \eqref{ba}
it is in fact quite useful, and if we multiply boths sides by $\alpha_n$ and 
then sum $n$ over the nonnegative integers we get
\begin{multline*}
\frac{q^{n(n+2)}}{(q;q)_{\infty}^3}
\sum_{n=0}^{\infty}\sum_{j=1}^{\infty}\sum_{r=0}^{\infty}
\sum_{l=0}^M \alpha_n (-1)^{M+j+r+1} 
q^{\binom{j+r}{2}+\binom{M-l}{2}+\binom{l}{2}+l(r+k+1)} \\
\times q^{-r(M-l)-n(2j-k)}
\frac{(1-q)(1-q^{(2j-k)(2n+1)})(1-q^{2r+k+1})}
{(1-q^{2n+1})(q;q)_{M-l}(q;q)_l}  \\[2mm]
=\begin{cases}
q^{m(m+k)}\beta_m & \text{$M=2m$} \\[2mm]
0 &\text{$M=2m+1$,} 
\end{cases}
\end{multline*}
where $(\alpha,\beta)$ is a Bailey pair relative to $q$.

Next we multiply both sides by $a^M$ and sum over $M$. If on the left
we interchange the sums over $M$ and $l$, shift $M\to M+l$ and then sum over
$l$ and $M$ using Euler's $q$-exponential sum \cite[Eq. (II.2)]{GR90}
\begin{equation*}
\sum_{n=0}^{\infty}\frac{(-1)^n a^n q^{\binom{n}{2}}}{(q;q)_n}=(a;q)_{\infty}
\end{equation*}
this yields
\begin{multline*}
\frac{(a;q)_{\infty}^2}{(q;q)_{\infty}^3}
\sum_{n=0}^{\infty}\sum_{j=1}^{\infty}\alpha_n q^{n(n-2j+k+2)}
\frac{(1-q)(1-q^{(2j-k)(2n+1)})}{(1-q^{2n+1})} \\
\times \sum_{r=0}^{\infty} \frac{(-1)^{j+r+1}
q^{\binom{j+r}{2}}(1-q^{2r+k+1})(aq^{-r};q)_r}{(a;q)_{r+k+1}}
=\sum_{n=0}^{\infty} a^{2n} q^{n(n+k)}\beta_n.
\end{multline*}
We have nearly arrived at \eqref{key}. All that is needed is the following
Bailey pair relative to $q$ due to Rogers \cite{Rogers17}:
\begin{equation*}
\alpha_n=(-1)^n q^{n(3n+1)/2}\frac{(1-q^{2n+1})}{(1-q)} 
\quad\text{and}\quad \beta_n=\frac{1}{(q;q)_n}.
\end{equation*}
Substituting this, interchanging the sum over $n$ and $j$ (with the
above choice for $\alpha_n$ this may indeed be done) and
using the triple product identity \eqref{tpi} gives \eqref{key}.

\begin{proof}[Proof of Lemma \eqref{lem1}]
Replacing $M$ by $2m+i$ where $i\in\{0,1\}$, and shifting $l\to l+m+i$
leads to
\begin{multline}\label{lquad}
\frac{q^{n(n+2)}}{(q;q)_{\infty}^3}
\sum_{j=1}^{\infty}\sum_{r=0}^{\infty}\sum_{l=-m-i}^m (-1)^{j+r+1}
q^{\binom{j+r}{2}+ir-n(2j-k)+l(l+2r+k+i+1)} \\
\times \frac{(1-q^{(2j-k)(2n+1)})(1-q^{2r+k+1})}
{(1-q^{2n+1})(q;q)_{m-l}(q;q)_{m+l+i}}
=\frac{\delta_{i,0}}{(q;q)_{m-n}(q;q)_{m+n+1}}.
\end{multline}
By the $q$-Chu--Vandermonde summation \cite[Eq. (II.6)]{GR90}
\begin{equation}\label{qCV}
\sum_{j=0}^n \frac{(a,q^{-n};q)_j \,q^j}{(q,c;q)_j}=
\frac{(c/a;q)_n}{(c;q)_n}\,a^n
\end{equation}
this follows from the simpler to prove identity
\begin{multline}\label{llin}
\frac{q^{2n(n+1)}}{(q;q)_{\infty}^3}
\sum_{j=1}^{\infty}\sum_{r=0}^{\infty}\sum_{l=-m-i}^m (-1)^{j+r+1}
q^{\binom{j+r}{2}+ir-n(2j-k)+l(2r+k+1)} \\
\times \frac{(1-q^{(2j-k)(2n+1)})(1-q^{2r+k+1})}
{(1-q^{2n+1})(q;q)_{m-l}(q;q)_{m+l+i}}
=\frac{\delta_{i,0}q^{m-n}}{(q;q)_{m-n}(q;q)_{m+n+1}}.
\end{multline}
Indeed, if we multiply both sides of \eqref{llin} by $q^{m(m+i)}/(q;q)_{M-m}$,
the resulting identity can be summed over $m$ by the $c=0$ instance of
\eqref{qCV} (after first replacing $m\to M-m$). On the right we of course
only need to do this sum when $i=0$. Replacing $M$ by $m$ then gives 
\eqref{lquad}. 
Those familiar with the concept of a Bailey chain 
\cite{Andrews84} will have recognized that the reduction of 
\eqref{lquad} to \eqref{llin} corresponds to a simplifying (i.e., backwards)
iteration along a Bailey chain relative to $q^i$.

Since \eqref{llin} is of the form \eqref{ab} with $a=q^i$ we can use 
\eqref{ba} to invert. Hence
\begin{multline}\label{inverted}
\frac{q^{2n(n+1)}}{(q;q)_{\infty}^3}
\sum_{j=1}^{\infty}\sum_{r=0}^{\infty}(-1)^{j+r+1}
q^{\binom{j+r}{2}+ir-n(2j-k)-(m+i)(2r+k+1)} \\
\times \bigl(1-q^{(2n+1)(2j-k)}\bigr)\bigl(1-q^{2r+k+1}\bigr)
\bigl(1+q^{(2m+i)(2r+k+1)}\bigr) \\
=\delta_{i,0}(1-q^{2m})(1-q^{2n+1})
\sum_{r=0}^m\frac{(-1)^{m-r}q^{\binom{m-r}{2}+r-n}(q;q)_{r+m-1}}
{(q;q)_{m-r}(q;q)_{r-n}(q;q)_{r+n+1}},
\end{multline}
with the convention that $(1-q^{2m})(q^m;q)_{r+m-1}=2$ for $m=r=0$
in accordance with $(1-q^{2m})(q;q)_{m-1}=(1+q^m)(q;q)_m$.
The sum over $r$ on the right may be carried out by the $q$-Chu--Vandermonde
sum \eqref{qCV}, leading to 
\begin{multline*}
\text{RHS}\eqref{inverted} \\ =
\delta_{i,0}\frac{(-1)^{m+n}(1-q^{2m})(1-q^{2\max\{n,-n-1\}+1})
(q;q)_{m+\max\{n,-n-1\}-1}}{(q;q)_{m-n}(q;q)_{m+n+1}
(q^2;q)_{\max\{n,-n-1\}-m}},
\end{multline*}
which is nonzero for $n=\pm m$ and $n=\pm m-1$ only.
If we also multiply \eqref{inverted} by $q^{i(k+1)/2}$ and 
note that on the right this may again be dropped, we obtain
\begin{multline}\label{del}
\frac{q^{2n(n+1)}}{(q;q)_{\infty}^3}
\sum_{j=1}^{\infty}\sum_{r=0}^{\infty}(-1)^{j+r+1} 
q^{\binom{j+r}{2}-n(2j-k)-(2m+i)(2r+k+1)/2} \\
\times \bigl(1-q^{(2n+1)(2j-k)}\bigr)\bigl(1-q^{2r+k+1}\bigr)
\bigl(1+q^{(2m+i)(2r+k+1)}\bigr) \\[1mm]
=\delta_{i,0}\bigl(\delta_{m,n}+\delta_{-m,n}-
\delta_{m-1,n}-\delta_{-m-1,n}\bigr).
\end{multline}
Since both sides are invariant under the substitution $m\to -m-i$ 
this must hold for all $m,n\in\Z$ and $i,k\in\{0,1\}$. 
 
Next we observe that \eqref{del} is a consequence of the stronger result
\begin{multline}\label{symm}
\frac{q^{\binom{n}{2}+(n-m)/2}}{(q;q)_{\infty}^3}
\sum_{j=1}^{\infty}\sum_{r=0}^{\infty}(-1)^{j+r+n}
q^{\binom{j+r}{2}-(n-1)(2j-k)/2-(m-1)(2r+k)/2} \\
\times \bigl(1-q^{n(2j-k)}\bigr)\bigl(1-q^{m(2r+k+1)}\bigr)
=\delta_{m,n}-\delta_{-m,n},
\end{multline}
for $m,n\in\Z$ and $k\in\{0,1\}$.
If we denote the above two identities by $\eqref{del}|_{m,n}$ and 
$\eqref{symm}|_{m,n}$ and note that $\delta_{2m+i\pm 1,2n+1}=
\delta_{i,0}\delta_{m-(1\mp 1)/2,n}$, then $\eqref{del}|_{m,n}=
\eqref{symm}|_{2m+i+1,2n+1}-\eqref{symm}|_{2m+i-1,2n+1}$.

Before proving \eqref{symm} let us point out that without loss of 
generality we may fix $k=0$. For, if we take \eqref{symm} with $k=1$, replace
$r\leftrightarrow j-1$, and multiply the result by $(-1)^{m-n}q^{(m^2-n^2)/2}$
we find $\eqref{symm}|_{m,n;k=1}=\eqref{symm}|_{n,m;k=0}$.
Equation \eqref{symm} for $k=0$ is a linear combination 
of yet another identity, given by
\begin{multline}\label{symm2}
\frac{q^{\binom{n}{2}}}{(q;q)_{\infty}^3}
\sum_{j=0}^{\infty}\sum_{r=0}^{\infty}(-1)^{j+r+n}
q^{\binom{j+r}{2}-(n-1)j-(m-1)r} \\
+\frac{q^{\binom{n}{2}}}{(q;q)_{\infty}^3}
\sum_{j=1}^{\infty}\sum_{r=0}^{\infty}(-1)^{j+r+n}
q^{\binom{j+r}{2}-(n-1)j-(m-1)r+2nj+m(2r+1)}
=\delta_{m,n}.
\end{multline}
Here it should be noted that the first sum over $j$ now includes the 
term $j=0$. It is easily seen that this extra term is cancelled out in the 
following linear combination, and that 
$\eqref{symm2}|_{m,n}-q^{-n}\eqref{symm2}|_{m,-n}=
q^{(m-n)/2}\eqref{symm}|_{m,n;k=0}$.
 
After this string of reductive steps we are finally in a position to
carry out a proof. Replacing $m$ by $m+n$ in \eqref{symm2} and changing
the summation variable $j\to n-j$ ($j\to j-n)$ in the first (second)
double sum gives
\begin{multline*}
\sum_{j=-\infty}^n\sum_{r=0}^{\infty}(-1)^{j+r}
q^{\binom{j-r}{2}-mr}
+\sum_{j=n+1}^{\infty}\sum_{r=0}^{\infty}(-1)^{j+r}
q^{\binom{j+r+1}{2}+m(r+1)} \\
=(q;q)_{\infty}^3 \delta_{m,0}.
\end{multline*}
In the second term on the left we rewrite the sum over $r$ using
\begin{equation}\label{zero1}
\sum_{r=0}^{\infty}(-1)^r q^{\binom{r+1}{2}+a(r+1)}
=\sum_{r=0}^{\infty}(-1)^r q^{\binom{r+1}{2}-ar}
\end{equation}
as follows from
\begin{equation}\label{zero2}
\sum_{r=-\infty}^{\infty}(-1)^r q^{\binom{r+1}{2}-ar}=0.
\end{equation}
(To prove \eqref{zero2} replace $r\to 2a-1-r$.) 
As a result we are left with
\begin{equation}\label{fin}
\sum_{j=-\infty}^{\infty}\sum_{r=0}^{\infty}(-1)^{j+r}
q^{\binom{j-r}{2}-mr}=(q;q)_{\infty}^3\delta_{m,0}.
\end{equation}
Using \eqref{zero1} on the sum over $r$ and negating $j$
yields $\eqref{fin}|_m=q^m \eqref{fin}|_{-m}$,
so that we may assume $m\leq 0$ when proving \eqref{fin}.
If $m<0$ the order of the sums may be interchanged. 
By \eqref{zero2} this completes the proof. 
If $m=0$ we need 
\begin{equation*}
\sum_{r=0}^{2j-1}(-1)^r q^{\binom{j-r}{2}}=0
\end{equation*}
for $j\geq 0$ (to prove this replace $r\to 2j-1-r$),
and Jacobi's identity \cite[\S 66, (5.)]{Jacobi29}
\begin{equation*}
\sum_{i=0}^{\infty}(-1)^i (2i+1) q^{\binom{i+1}{2}}=(q;q)_{\infty}^3.
\end{equation*}
Equipped with these the rest is easy;
\begin{align*}
\sum_{j=-\infty}^{\infty}\sum_{r=0}^{\infty}(-1)^{j+r}q^{\binom{j-r}{2}}
&=\Bigl\{\sum_{j=0}^{\infty}\sum_{r=2j}^{\infty}+
\sum_{j=-\infty}^{-1}\sum_{r=0}^{\infty}\Bigr\}
(-1)^{j+r} q^{\binom{j-r}{2}} \\
&=\Bigl\{\sum_{j=0}^{\infty}
+\sum_{j=1}^{\infty}\Bigr\} \sum_{r=j}^{\infty}
(-1)^r q^{\binom{r+1}{2}}\\
&=\sum_{r=0}^{\infty}(-1)^r (2r+1) q^{\binom{r+1}{2}}
=(q;q)_{\infty}^3. \qquad \qed
\end{align*}
\renewcommand{\qed}{}
\end{proof}

\section{Discussion}\label{secd}
\subsection{Eqs~\eqref{RRP1} and \eqref{RRP2} versus \eqref{AJ}}
The proof of Theorem~\ref{thm} as given in the previous section is very
lengthy and complicated, and, as a result, not very illuminating.
Here we briefly discuss the proofs of \eqref{AJ} as found by Jordan and
Andrews as we hold some hope that at least one of these may be generalized to 
also prove \eqref{RRP1} and \eqref{RRP2}.

Perhaps simplest is Jordan's proof~\cite{Jordan84}.
Denoting the right side of \eqref{AJ} by $f_{n;k}$ and the summand on 
the right of \eqref{AJ} by $f_{n,r;k}$, it is not difficult to show that
the functional equation 
\begin{equation*}
(1-x q^{n+2})K_{n+1}(x)=K_n(x)-x^2 q^{n+4}K_n(x q^2)
\end{equation*}
satisfied by the Szeg\"o polynomials implies the recurrence
\begin{equation}\label{frec}
\sum_{r=m+1}^n (f_{n,r;k}-f_{n-1,r-1;k})=
-\frac{1-q^{m-n}}{1-q^{n+k}}f_{n,m;k}.
\end{equation}
By the $m=0$ instance hereof it is found that
\begin{align*}
f_{n;k}-f_{n-1;k}&=f_{n,0;k}+\sum_{r=1}^n(f_{n,r;k}-f_{n-1,r-1;k}) \\
&=q^{-n}\frac{1-q^{2n+k}}{1-q^{n+k}}f_{n,0;k} \\[2mm]
&=(-1)^n q^{n(5n+2k+1)/2}+(-1)^{n+k}q^{(n+k)(5(n+k)-2k-1)/2},
\end{align*}
from which \eqref{AJ} follows by induction.
Unfortunately, at present we have been unable to find an analogue of \eqref{frec}
for the summands on the right of Theorem~\ref{thm}.

Andrews' proof of \eqref{AJ} relies on the following multiple series 
generalization of Watson's $q$-Whipple transform 
\cite[Thm. 4; $k=3$]{Andrews75}:
\begin{multline}\label{k3}
{_{10}}W_9(a;b,c,d,e,f,g,q^{-n};q,a^3 q^{n+3}/bcdefg;q,q) \\
=\frac{(aq,aq/fg;q)_n}{(aq/f,aq/g;q)_n}
\sum_{j=0}^n \sum_{k=0}^{n-j}
\frac{(aq/bc,d,e;q)_j}{(q,aq/b,aq/c;q)_j}
\frac{(aq/de;q)_k}{(q;q)_k} \\ \times
\frac{(f,g,q^{-n};q)_{j+k}}{(aq/d,aq/e,fgq^{-n}/a)_{j+k}}
\Bigl(\frac{aq^2}{de}\Bigr)^j q^k.
\end{multline}
Here and in the following we employ standard notation for basic 
hypergeometric series, see e.g.,~\cite{GR90}.
Taking $b=aq^{n+1}$ and letting $c,d,e,f,g$ tend to infinity 
yields \eqref{AJ} with $k=0$ if $a=1$ ($k=1$ if $a=q$).

Now if we apply Sears' $_{4}\phi_3$ transformation \cite[Eq. (III.15)]{GR90}
to Watson's  Watson's $q$-Whipple transform \cite[Eq. (III.18)]{GR90}
we readily obtain
\begin{multline}\label{WS}
{_8W_7}(a;b,c,d,e,q^{-n};q,a^2 q^{n+2}/bcde;q,q) \\
=\frac{(aq,b,a^2 q^2/bcde;q)_n}{(aq/c,aq/d,aq/e;q)_n}
{_4\phi_3}\Bigl[\genfrac{}{}{0pt}{}
{aq/bc,aq/bd,aq/be,q^{-n}}{aq/b,a^2q^2/bcde,q^{1-n}/b};q,q\Bigr].
\end{multline}
Taking $b=aq^{n+1}$ and letting $c,d,e,f,g$ this simplifies to
\begin{equation}\label{WSlim}
\sum_{j=0}^n \frac{1-aq^{2j}}{1-a}
\frac{(a;q)_j(-1)^j q^{3\binom{j}{2}}(aq)^j}
{(q;q)_j}=(aq;q)_{2n}\sum_{j=0}^n \frac{q^j}{(q,q^{-2n}/a;q)_j}.
\end{equation}
Chosing $a=q^k$ for $k=0,1$ and making the variable change $j\to n-r$ on the right
we obtain the following polynomial analogue of Euler's identity
\begin{equation}\label{Eulerpol}
\sum_{j=-n-k}^n(-1)^j q^{j(3j+1)/2}=\sum_{r=0}^n (-1)^{n-r}
q^{(n-r)(3n+r+2k+3)/2}\frac{(q;q)_{n+r+k}}{(q;q)_{n-r}}.
\end{equation}
(Incidentally, this identity is very similar and can easily be transformed
into a polynomial version of Euler's identity due to Shanks \cite{Shanks51}.)
Given the similarity between \eqref{Eulerpol} and the identities \eqref{RRP1} and 
\eqref{RRP2}, and given Andrews' proof of \eqref{AJ} by means
of \eqref{k3} it seems very natural to ask for a proof of Theorem~\ref{thm}
by means of a multiple series generalization of the transformation~\eqref{WS}.
If we take \eqref{RRP1} with $k=0$ and \eqref{RRP2} with $k=1$ 
and replace $r\to n-j$ in the sums on the right we find that the resulting
identities are the $a=1$ and $a=q$ instances of
\begin{multline*}
\sum_{j=0}^n \frac{1-aq^{2j}}{1-a}
\frac{(a;q)_j(-1)^j q^{5\binom{j}{2}}(aq)^{2j}}{(q;q)_j} \\
=(aq;q)_{2n}\sum_{j=0}^n\sum_{k=0}^{n-j}
\frac{(-1)^k a^{j+k} q^{j^2-\binom{2j+k}{2}+(2n+1)(j+k)}
(q^{-2n-1};q)_{2j+2k}}
{(q,q^{-2n}/a;q)_j(q;q)_k(q^{-2n-1};q)_{2j+k}}
\end{multline*}
This is to be compared with \eqref{WSlim}.
Despite numerous attempts we failed to extend this to a multiple
series transformation similar to \eqref{k3} and generalizing \eqref{WS}.
Of course one can try to prove the above by equating coefficients of $a^m$,
but the resulting identity 
\begin{multline*}
\sum_{j=0}^n \frac{(-1)^j q^{\binom{j}{2}+j(4j-2m+1)}}{(q;q)_j}
\biggl(\qbin{j}{m-2j}+q^{4j-m+1}\qbin{j}{m-2j-1}\biggr) \\
=\sum_{r=0}^n\sum_{s=0}^{n-r}
\frac{(-1)^{r+s}q^{\binom{r+s}{2}+r(4n-2m+4)+s(s+r-m+1)}}{(q;q)_r}
\\ \times
\qbin{2n+1-r}{m-2r-s}\qbin{2n-2r-s+1}{s}
\end{multline*}
for $n\geq 0$ and $0\leq m\leq 3n+1$ is not particularly simple
(and would only prove half of Theorem~\ref{thm}).

\subsection{Some combinatorics related to Theorem~\ref{thm}}
In order to discuss some of the combinatorics of Theorem~\ref{thm}
we need to review several standard results from partition 
theory \cite{Andrews76}.

Let $\lambda=(\lambda_1,\lambda_2,\dots,\lambda_r)$ be a partition, 
defined as a weakly decreasing sequence of positive integers $\lambda_j$ 
(the parts of $\lambda$). The weight $|\lambda|$ of $\lambda$
is given by the sum of its parts. We say that $\lambda$ is a partition of
$l$ if $|\lambda|=l$. 
The Ferrers graph of $\lambda$ is the graph obtained
by drawing $r$ left-aligned rows of dots with the $j$th row containing
$\lambda_j$ dots. 
The conjugate $\lambda'$ of $\lambda$ is obtained by transposing its
Ferrers graph. 
The number $d(\lambda)$ is the number of rows in the
maximal square of dots of the Ferrers graph of $\lambda$.
An alternative way to represent a partition $\lambda$
is as a two-rowed matrix of $d(\lambda)=d$ columns
$\bigl(\begin{smallmatrix} t_1 t_2 \dots t_d \\ b_1 b_2 \dots b_d
\end{smallmatrix}\bigr)$, where $t_j=\lambda_j-j$ and $b_j=\lambda_j'-j$, 
so that, in particular, $t_j>t_{j+1}$ and $b_j>b_{j+1}$. Conversely, any 
such matrix (also called Frobenius symbol) corresponds to the unique
partition $\lambda$ by $\lambda_j=t_j+b_j+1$.
We will in the following identify the 
standard and Frobenius notations for partitions. Note that 
$|\lambda|=d+\sum_{j=1}^d(t_j+b_j)$.
The rank of a partition $\lambda$ is defined as its largest part minus 
its number of parts, i.e., as $\lambda_1-\lambda'_1=t_1-b_1$.
More generally, the $i$th successive rank of $\lambda$ is given by 
$t_i-b_i$, and $r(\lambda)=(t_1-b_1,t_2-b_2,\dots,t_d-b_d)$ denotes 
the sequence of successive ranks of $\lambda$.
For example, if $\lambda=(7,7,5,3,3,1,1,1)$, then $|\lambda|=28$, 
$\lambda'=(8,5,5,3,3,2,2)$, $d(\lambda)=3$,
$\lambda=\bigl(\begin{smallmatrix} 6 5 2 \\ 7 3 2 \end{smallmatrix}\bigr)$, 
$\lambda'=\bigl(\begin{smallmatrix} 7 3 2 \\ 6 5 2 \end{smallmatrix}\bigr)$,
and $r(\lambda)=(-1,2,0)$.

Now let $b_2(l,n)$ denote the set of all partitions of $l$, with largest 
part at most $n-2$ and difference between parts at least $2$,
and let $B_2(l,n)$ be its cardinality. Then
$e_n=\sum_{l\geq 0}B_2(l,n) q^l$.
Given a partition $\lambda\in b_2(l,n)$ with exactly $r$ parts, 
one can form a new partition $\mu$ as follows \cite[\S 9.3]{Andrews76}:
$\mu=\bigl(\begin{smallmatrix} s_1,\dots,s_r \\ c_1,\dots,c_r 
\end{smallmatrix}\bigr)$,
where $s_j=\lfloor \lambda_j/2 \rfloor$ and 
$c_j=\lfloor (\lambda_j-1)/2 \rfloor$. 
Because of the difference-$2$ condition one indeed has $s_j>s_{j+1}$ and
$c_j>c_{j+1}$. Since (for $n\in\Z$) 
$\lfloor n/2 \rfloor +\lfloor (n-1)/2 \rfloor=n-1$ one 
finds that $|\mu|=r+\sum_{j=1}^r (s_j+c_j)=|\lambda|=l$.
Furthermore, the restriction that $\lambda_j-\lambda_{j+1}\geq 2$ translates
into the fact that the successive ranks of $\mu$ must take the values $0$ and
$1$ only. Finally the restriction that $\lambda_1\leq n-2$ implies that
$s_1+c_1+1\leq n-2$. Since $s_1-c_1\in\{0,1\}$ this is equivalent to
requiring that $\mu_1\leq\lfloor n/2\rfloor$ and 
$\mu'_1\leq\lfloor (n-1)/2\rfloor$.
If we denote the set of all partitions of $l$ with successive ranks in 
$\{0,1\}$, largest part not exceeding $\lfloor n/2 \rfloor$ and number 
of parts not exceeding $\lfloor (n-1)/2\rfloor$ by $q_2(l,n)$ 
(with cardinality $Q_2(l,n)$) then clearly each partition 
$\mu\in q_2(l,n)$ can also be mapped back onto a partition in $b_2(l,n)$.
Specifically, if $\mu\in q_2(l,n)$ has Frobenius symbol
$\bigl(\begin{smallmatrix} s_1,\dots,s_r \\ c_1,\dots,c_r 
\end{smallmatrix}\bigr)$,
then $\lambda=(s_1+c_1+1,\dots,s_r+c_r+1)\in  b_2(l,n)$, since
$\lambda_j-\lambda_{j+1}=s_j+c_j-s_{j+1}-c_{j+1}\geq 2$ and 
$\lambda_1=s_1+c_1+1=\mu_1+\mu'_1-1\leq n-2$.
Hence $Q_2(l,n)=B_2(l,n)$ and $e_n=\sum_{l\geq 0}Q_2(l,n) q^l$.
For example, $\cup_{l\geq 0}q_2(l,n)=
\{\emptyset,(1),(2),(2,1),(3,1),(2,2),(3,2),(3,3)\}$ so that
$e_6=1+q+q^2+q^3+2q^4+q^5+q^6$.

The above discussion can be repeated for the Schur polynomial $d_n$ and
we define $b_1(l,n)$ as the subset of $b_2(l,n)$ obtained by
removing all partitions which have a part equal to $1$. Hence
$d_n=\sum_{l\geq 0}B_1(l,n)q^l$.
If we also define $q_1(m,n)$ as the set of partitions with successive ranks in
$\{1,2\}$, largest part not exceeding $\lfloor (n+1)/2 \rfloor$ and number of 
parts not exceeding $\lfloor (n-2)/2\rfloor$,
then it is not hard to show that $Q_1(l,n)=B_1(l,n)$ so that
$d_n=\sum_{l\geq 0} Q_1(l,n) q^l$.
For example, $\cup_{l\geq 0} q_1(l,n)=\{\emptyset,(2),(3),(3,1),(3,3)\}$ 
so that $d_6=1+q^2+q^3+q^4+q^6$.

So far, we have given a combinatorial interpretation of the
Schur polynomials $e_n$ and $d_n$ in terms of partitions with restrictions
on their size and successive ranks. Next we will discuss the combinatorial
interpretation of the partial theta sum 
$\sum_{j=-n+k}^n (-1)^j q^{j(5j-2i+5)/2}$
in terms of successive ranks.

First we recall some further known properties of $Q_i(l,n)$ 
\cite{Andrews72,Andrews76}.
Let $\lambda$ be a partition and $r(\lambda)$ its sequence of successive ranks.
The length of the largest subsequence $r'$ of $r(\lambda)$ such that 
the odd (even) elements of  $r'$ are at least $4-i$ and the even (odd) elements
of $r'$ are at most $1-i$, is called the $(2,i)$-positive ($(2,i)$-negative)
oscillation of $\lambda$.
The number of partitions of $l$ that have at most $b$ parts, largest part
not exceeding $a$ and $(2,i)$-positive ($(2,i)$-negative) oscillation at least
$j$ is denoted by $p_i(a,b;j;l)$ ($m_i(a,b;j;l)$).
By inclusion-exclusion arguments it then follows that
\begin{equation*}
Q_i(l,n)=\sum_{j=0}^{\infty}(-1)^j p_i(\bar{a},\bar{b};j,l)
+\sum_{j=1}^{\infty}(-1)^j m_i(\bar{a},\bar{b};j,l),
\end{equation*}
with $\bar{a}=\bar{a}(n,i)=\lfloor (n-i+2)/2\rfloor$ 
and $\bar{b}=\bar{b}(n,i)=\lfloor (n+i-3)/2\rfloor$.
Furthermore,
\begin{equation*}
q^{j(5j-2i+5)/2}\qbin{n-1}{\lfloor \frac{n+i-5j-3}{2}\rfloor}=
\begin{cases}
\sum_{l=0}^{\infty} p_i(\bar{a},\bar{b};-j,l)q^l 
&\text{$j\leq 0$, $j$ even} \\[1mm]
\sum_{l=0}^{\infty} m_i(\bar{a},\bar{b};-j,l)q^l 
&\text{$j\leq 0$, $j$ odd} \\[1mm]
\sum_{l=0}^{\infty} p_i(\bar{a},\bar{b};j,l)q^l 
&\text{$j\geq 0$, $j$ odd} \\[1mm]
\sum_{l=0}^{\infty} m_i(\bar{a},\bar{b};j,l)q^l 
&\text{$j\geq 0$, $j$ even,}
\end{cases}
\end{equation*}
from which \eqref{edbos} immediately follows.
But now it is also clear what our partial theta sums represent.
If we denote by $\lambda^{\pm}_{i,j}$ the
(unique) partition of minimal weight that has a positive/negative
$(2,i)$-oscillation $j$ and by $M_i$ the set of all such 
minimal partitions, i.e.,
$M_i=\{\lambda^{\sigma}_{i,j}\}_{j\geq 0;\sigma\in\{0,1\}}$, then
(for $k\in\{0,1\}$ and $i\in\{1,2\}$)
\begin{equation*}
\sum_{j=-n-k}^n (-1)^j q^{j(5j-2i+5)/2}=
\sum_{\substack{\lambda\in M_i \\
\lambda_1\leq \lfloor (5n+2ki-2i+5)/2\rfloor \\
\lambda'_1\leq \lfloor (5n+2ki)/2\rfloor}} 
(-1)^{d(\lambda)} q^{|\lambda|}.
\end{equation*}
One can in fact easily find the partition $\lambda^{\pm}_{i,j}$.
For example, using the Frobenius notation it follows immediately that
for $j$ even
\begin{equation*}
\lambda^{+}_{i,j}=
\begin{pmatrix} 5j/2-1,&5j/2-5,&\dots,& 9,& 5, &4, & 0 \\ 
5j/2+i-5,&5j/2+i-6,& \dots,& i+5,& i+4,& i,& i-1\end{pmatrix}.
\end{equation*}
When converted into standard notation this gives
\begin{equation*}
\lambda^{+}_{i,j}=
(5j/2,(5j/2-3)^2,\dots,(j+6)^2,(j+3)^2,j^i,(j-2)^3,\dots,4^3,2^3)
\end{equation*}
where $p^f$ stands for $f$ parts of size $p$. 
Calculating the weight of this partition gives
\begin{equation*}
|\lambda^{+}_{i,j}|=
5j/2+2\sum_{k=1}^{j/2-1}(j+3k)+ij+3\sum_{k=1}^{j/2-1}(2k)=j(5j+2i-5)/2,
\quad \text{$j$ even}
\end{equation*}
as it should. Similarly one can use the Frobenius notation to find
\begin{equation*}
\lambda^{-}_{i,j}=
((5j-3)/2)^2,\dots,(j+6)^2,(j+3)^2,j^i,(j-2)^3,\dots,3^3,1^3),
\end{equation*}
$(\lambda^{+}_{i,j})'=\lambda^{-}_{5-i,j}$ both for $j$ odd,
and $(\lambda^{-}_{i,j})'=\lambda^{+}_{5-i,j}$ for $j$ even, and thus
$|\lambda^{\pm}_{i,j}|=j(5j\mp 2i \pm 5)/2$ for odd $j$ and
$|\lambda^{-}_{i,j}|=j(5j-2i+5)/2$ for even $j$.

Summarizing, we have the following remarkable situation.
The Schur polynomials, which are the generating functions of certain
size and successive rank restricted partitions, can be expressed as 
an alternating sum over the generating functions of partitions with 
certain restrictions on their $(2,i)$-oscillations.
This well-known fact \cite{Andrews72,Andrews76} provides a combinatorial
explanation of Schur's result \eqref{edbos}. 
But now we see that according to the Theorem~\ref{thm} there is another 
side to the coin;
the alternating sum over the generating function of a very special subset 
of partitions with certain restrictions on their $(2,i)$-oscillations
can in its turn be expressed as a weighted sum over Schur polynomials.
However, by no means is this an example of a trivial (or nontrivial but
known) inversion result. Indeed, naively one might think that if 
we substitute 
\eqref{edbos} in Theorem~\ref{thm} to get
\begin{multline*}
\sum_{j=-n-k}^n(-1)^j q^{j(5j-2i+5)/2}=
\sum_{j=-\infty}^{\infty}(-1)^j q^{j(5j-2i+5)/2} \\
\times \sum_{r=0}^n (-1)^{n-r}
q^{(n-r)(5n+3r+4k+5)/2}\frac{(q;q)_{n+r+k}}{(q;q)_{n-r}} 
\qbins{2r+k+1}{\lfloor{(2r+i+k-5j-1)/2}\rfloor},
\end{multline*}
that this is just a consequence of the second line being 
$\chi(-n-k\leq j\leq n)$ with $\chi(\text{true})=1$ and 
$\chi(\text{false})=0$. However, it is readily checked that this is only 
correct when $n=k=0$.
It thus seems an extremely challenging problem to find a combinatorial
proof of Theorem~\ref{thm}, especially since our analytic
proof provides so little insight as to why this theorem is true.

To conclude we remark that the previous discussion
has a representation theoretic counterpart. As is well-known, 
\begin{equation*}
\frac{1}{(q;q)_{\infty}}
\sum_{j=-\infty}^{\infty}(-1)^j q^{j(5j-2i+5)/2}
\end{equation*}
is the (normalized) character of the $c=-22/5$ Virasoro
algebra corresponding to the highest weight vector $v_{h_i}$ of
weight $h_i=(1-i)/5$. 
According to the Feigin--Fuchs construction \cite{FF84} the above character
can be constructed from the Verma module $V(c,h_i)$ by eliminating submodules
generated by singular or null vectors. Because of the embedded
structure of these submodules this leads to an inclusion-exclusion
type of sum. Specifically, the character corresponding to
the submodule $V(c,h_i')$ with singular vector of weight $h_i'$
is given by  $q^{h_i'-h_i}/(q;q)_{\infty}$, with the
set of weights of singular vectors (including $v_{h_i}$)
given by $h_i'=h_i+j(j-2i+5)/2$ for $j\in\Z$.
Therefore, if we denote by $V_s(c,h_i;N)$ the set comprising of the 
$N$ singular vectors
of $V(c,h_i)$ of smallest weight, and if we denote by $d(v)+1$ the number of
(sub)modules $V(c,h_i')$ that contain the singular vector $v$ 
(so that $d(v)=0$ iff $v=v_{h_i}$), then
\begin{equation*}
\sum_{j=-n-k}^n (-1)^j q^{j(5j-2i+5)/2}=
\sum_{v\in V_s(c,h_i;2n+k+1)} (-1)^{d(v)} q^{|v|-h_i},
\end{equation*}
where $|v|$ is the weight of $v$, $c=-22/5$ and $h_i=(1-i)/5$.
Again it is a challenge to explain Theorem~\ref{thm} from the above
representation theoretic point of view.

\subsection*{Note added}
Robin Chapman has informed me that he has found a combinatorial
proof of Theorem~\ref{thm} using Schur's involution.

\bibliographystyle{amsplain}

\end{document}